\documentclass[10pt]{article}
\usepackage{amsmath, amsthm, amssymb, amsfonts, latexsym,amscd}
\usepackage{authblk}
\usepackage{tikz}
\usepackage{mathtools}
\usetikzlibrary{shadings,patterns}
\newtheorem{theorem}{Theorem}[section]
\newtheorem{corollary}[theorem]{Corollary}
\newtheorem{lemma}[theorem]{Lemma}
\newtheorem{proposition}[theorem]{Proposition}
\newtheorem{definition}[theorem]{Definition}
\newtheorem{remark}[theorem]{Remark}
\newtheorem{example}[theorem]{Example}
\newtheorem{question}[theorem]{Question}

\newtheorem{conjecture}[theorem]{Conjecture}
\begin{document}

\title{\bf Some new central parts of connected  graphs}

\author[1,2]{Dinesh Pandey \thanks{Corresponding Author: dinesh.pandey@niser.ac.in, Supported by UGC Fellowship scheme (Sr. No. 2061641145), Government of India } }
\author[1,2]{Kamal Lochan Patra \thanks{klpatra@niser.ac.in}}
\affil[1]{School of Mathematical Sciences,
National Institute of Science Education and Research (NISER), Bhubaneswar,
P.O.- Jatni, District- Khurda, Odisha - 752050, India 
}

\affil[2]{Homi Bhabha National Institute (HBNI),
Training School Complex, Anushakti Nagar,
Mumbai - 400094, India}
\date{}
\maketitle
\begin{abstract}
The center, median and the security center are three central parts defined for any connected graph whereas the characteristic set, subtree core and core vertices are three central parts defined for trees only. We extend the concept of the characteristic set, subtree core and core vertices to general connected graphs and call them the characteristic center, subgraph core and core vertices, respectively. 

We show by examples that  in a connected graph all the above six central parts can be different and also prove that for a connected vertex transitive graph each of the six central parts is the whole vertex set. Further it is shown that given any graph $G$, there  exists a connected supergraph $G_{ch}$ of $G$ with the whole vertex set of $G$ as the characteristic center. Associated with the subgraph core and core vertices, we leave some unanswered question related to the graph centrality. \\


\end{abstract}


\section{Introduction}
Throughout this paper, graphs are simple, finite, connected and undirected. Let $G$ be a graph with vertex set $V(G)$ and edge set $E(G)$. A vertex $v\in V(G)$ is called a {\it cut vertex} of $G$ if the graph $G-v$ is disconnected. If a graph has no cut vertices, it is called a two connected graph. A maximal two connected subgraph of $G$ is called a {\it block} of $G$. For $u,v \in V(G),$ the {\it distance} between them is defined as the number of edges in a shortest path joining $u$ and $v$, and we denote it by $d(u,v)$. 

In $1869,$ Jordan first introduced the notion of centrality in connected graphs. In \cite{Jor}, he first defined two central parts for trees, popularly known as the {\it center} and the {\it centroid}. The center has a natural extension to connected graphs while the median and the security center are considered as two different generalizations of the centroid in connected graphs. Including centroid there are many other central parts which are defined for trees only. Some of these are the subtree core, the core vertices, the characteristic set, the telephone center, the weight balance center, the processing center, the n-th power center of gravity, the k-nucleus etc.  Interested people may see the survey paper by Reid \cite{Reid} for many other central parts defined for a tree. 

Central parts of graphs have practical importance in finding an optimal location for establishing a facility in a network. Many facility location problems can be solved by translating it into a graph theoretic problem. Most of the time the solutions of such problems are given by a central part of the associated graph. The associated graph may not be a tree always. So, it is worth to extend some of the definitions of the central parts of trees  to connected graphs.  

We denote the automorphism group of a graph $G$ by $Aut(G)$. A graph $G$ is called a {\it vertex transitive graph} if for any two vertices $u$ and $v$ of $G$, there exists $\sigma \in Aut (G)$ such that $\sigma(u)=v$. We denote $\sigma(u)$ by $u^{\sigma}$. It is to note that a vertex transitive graph is always regular.

\begin{lemma}[\cite{God}, Lemma 1.3.2]\label{distance preserve}
Let $u,v\in V(G)$ and $\sigma\in Aut(G)$. Then $d(u,v)=d(u^{\sigma},v^{\sigma})$.
\end{lemma}

The characteristic set, subtree core and  core vertices  are  three central parts defined for trees only. This paper is a study about the generalization of these three central parts to connected graphs. In Section \ref{Prim}, we see some basic properties of the central parts: center, median and security center. In Section \ref{New centers}, we define three new central parts of a connected graph which are generalizations of the characteristic set, subtree core and  core vertices. We explore their properties analogous to  the center,  median and security center. While exploring the properties of the central parts of connected graphs, we talk about many unanswered questions related to the newly defined central parts.

\section{Preliminaries}\label{Prim}
In this section, we recall the definitions of the central parts: center, median and security center and observe some properties which are common to them. 
\subsection{Center of a graph}\label{center}

The {\it eccentricity} $e(v)$ of a vertex $v$ of $G$ is defined as $e(v)=\max\{d(v,u): u\in V(G)\}$. A vertex of minimum eccentricity is called a {\it central vertex} of $G$  and the set of all central vertices is called the {\it center} of $G$. We denote the center of $G$ by $C(G)$.  For any tree $T$, the following is known regarding $C(T).$

\begin{proposition} [\cite{H}, Theorem 4.2] \label{Tree center}
The center of a tree  consists of either one vertex or two adjacent vertices.
\end{proposition}
 
The above result is generalized by  Harary and Norman in \cite{Hn}. 
 
\begin{proposition}[\cite{Hn}, Lemma 1]
 The center of a  graph $G$ is contained in a block of $G$.
\end{proposition}

A graph $G$  is called a {\it self-centered} graph if  $C(G)=V(G)$. By Proposition \ref{Tree center} it is clear that $K_1$ and $K_2$ are the only self-centered trees.  It will be nice to characterize all self-centered connected graphs. We find a class of graphs which are self-centered. For the complete graph $K_n$, the eccentricity of every vertex is $1$, so  $C(K_n)=V(K_n)$. Similarly the eccentricity of every vertex of the cycle $C_n$ is $\lfloor\frac{n}{2}\rfloor$ and so $C(C_n)=V(C_n)$.  Both $K_n$ and $C_n$ are vertex transitive graphs. Interestingly, it is true that any connected vertex transitive graph is self-centered.

\begin{theorem}\label{C_trans}
Let $G$ be a connected vertex transitive graph. Then $C(G)=V(G)$.
\end{theorem}

\begin{proof}
Let $G$ be a connected vertex transitive graph and let $u,v \in V(G)$.  Suppose $e(u)=d(u,u')$ and $e(v)=d(v,v')$. Since $G$ is vertex transitive,  there exists $\sigma \in Aut(G)$ such that $v^{\sigma}=u$.  By Lemma \ref{distance preserve}, we have $e(u)\geq d(u,v'^\sigma)=d(v^\sigma, v'^\sigma)=d(v,v')=e(v)$. Similarly, it  can be showed that $e(v)\geq e(u)$. So $e(u)=e(v)$ and hence $C(G)=V(G)$.
\end{proof}

The converse of  Theorem \ref{C_trans} is not true. For example, consider  the complete bipartite graph $K_{m,n}$ with $m> n\geq 2$. For any vertex $v\in V(K_{m,n})$, $e(v)=2$ and so $C(K_{m,n})=V(K_{m,n})$. But $K_{m,n}$  is not vertex transitive as it is not regular. 

We call the subgraph induced by the central vertices of $G$ as the {\it center subgraph} of $G$. The next property of the center is due to Buckley, Miller and Slater.

\begin{theorem} [\cite{Buc}]
For any graph $G$ (may be disconnected), there exists a connected graph $G'$ such that the center subgraph of $G'$ is isomorphic to $G$.
\end{theorem}

\subsection{Median of a graph}\label{median}

The median of a connected graph was first defined by Ore (see \cite{O}, page $30$). It was Zelinka who  considered it as a generalization of the centroid for connected graphs and observed it  as a central part of a connected graph. We start this subsection with the definition of the centroid of a tree.

Let $T$ be a tree. For $v\in V(T)$, a $\it{branch}$ at $v$ is a maximal subtree of $T$ containing $v$ as a pendant vertex. The {\it weight} of $v$ is the maximal number of edges in any branch of $T$ at $v.$  A vertex of minimal weight is called a centroid vertex of $T$ and the set of all centroid vertices is called the {\it centroid} of $T.$ It is known that the centroid of a tree consists of either one vertex or two adjacent vertices.

The distance $D(v)$ of a vertex $v$ in $G$ is defined as $D(v)=\sum_{u\in V(G)}d(v,u)$. A vertex of minimum distance is called a {\it median vertex} of $G$ and the set of all median vertices is called the {\it median} of $G$. We denote the median of $G$ by $M(G)$. In \cite{Z}, Zelinka proved that, while restricted to trees, the median coincides with the centroid. So, the median is considered as an extension of the centroid in connected graphs. Due to Zelinka, we have  the following property of the median analogous to Proposition \ref{Tree center}. 

\begin{proposition}[\cite{Z}, Theorem 2 ]\label{prop:med1}
The median of a tree consists of either one vertex or two adjacent vertices. 
\end{proposition}

The above result is generalized by Truszczy\'nski in \cite{T}. 

\begin{proposition}[\cite{T}, Theorem 3]
The median of a connected graph $G$ is contained in a block of $G$.
\end{proposition}

Next we prove the result corresponding to Theorem \ref{C_trans} for median. 
\begin{theorem}\label{Cd_trans}
Let $G$ be a vertex transitive graph. Then $M(G)=V(G)$.
\end{theorem}

\begin{proof} 
Let $u,v\in V(G)$. It is sufficient to prove that  $D(u)=D(v)$. Since $G$ is vertex transitive, there exists $\sigma \in Aut(G)$ such that $u^{\sigma}=v$.  By Lemma \ref{distance preserve}, we have  $D(u)= \underset{x\in V(G)}\sum {d(u,x)} = \underset{x\in V(G)}\sum {d(u^{\sigma},x^{\sigma}})$. So, 
$$D(u)=  \underset{x\in V(G)}\sum {d(u^{\sigma},x^{\sigma})} = \underset{x\in V(G)}\sum {d(v,x^{\sigma})}=\underset{y\in V(G)}\sum {d(v,y)}= D(v).$$
\end{proof}

Is the converse of Theorem \ref{Cd_trans} true?  To be more precise, we have the following question.
\begin{question}
 Let $G$ be a connected graph with $M(G)=V(G)$. Is $G$  vertex transitive? 
\end{question}

 The subgraph induced by the median vertices of $G$ is said to be the {\it median subgraph} of $G$. Due to Slater we have two different supergraphs whose median subgraphs are isomorphic to a given graph.
 
 \begin{proposition}[\cite{Sla2}, Theorem 2]\label{median supergraph}
  Let $G$ be a graph (may be disconnected). Then there exists a connected graph $G'$ such that the median subgraph of $G'$ is isomorphic to $G$.
  \end{proposition}
  
   The construction in Proposition \ref{median supergraph} contains  a large number of vertices. In \cite{Smart}, Smart and Slater gave another construction with less number of vertices  for graphs with no isolated vertices.
  
\begin{proposition}[ \cite{Smart}, Theorem 7]
 If $G$ has no isolated vertices, then there exists a connected graph $G'$ with $|V(G')|=2|V(G)|$ such that the median subgraph of $G'$ is isomorphic to $G$.
 \end{proposition}
 
 \subsection{Security center of a graph}\label{Scenter} 
 
The security center of a connected graph was introduced as another generalization of the  centroid in \cite{Sla}. For distinct vertices $u$ and $v$ of $G$, let $V_{uv}=\{x\in V(G):d(x,u)<d(x,v)\}$ and let $g(u,v)=|V_{uv}|-|V_{vu}|$. The {\it security number} $s(u)$ of $u$ is given by $s(u)=\min\{g(u,v):v\in V(G)-u\}$. The {\it security  center} of $G$, denoted by $\mathbb{S}(G)$ is the set of vertices $x$ for which $s(x)$ is maximum. 

The weight of a vertex $v$ in a tree $T$ is independent of the distance $d(u,v)$ for $u$ lying in a branch at $v$ rather it depends upon the number of vertices present in a branch at $v$. By this observation, Slater (\cite{Sla}) felt that the security center would be a better generalization of the centroid than the median. He proved that the security center of a tree coincides with its centroid. Due to Slater we have the following. 
 
\begin{proposition}[\label{S_center=centroid}\cite{Sla}, Corollary 1a]
The security center of a tree consists of either one vertex or two adjacent vertices .
\end{proposition}

Slater has shown by an example (\cite{Sla}, Figure 4 ) that, in a graph $G$ the  security center can be different from the median but both lie in the same block of $G$. Thus, regarding the position of the security center in a graph, we have the following result.

\begin{proposition}[\label{S_centerblock} \cite{Smart}, Theorem 6]
The security center of a graph $G$ is contained in a block of $G$.
\end{proposition}

Next we show the centrality property of vertex transitive graphs for the security center.
\begin{lemma}\label{g(u,v)}
Let $u$ and $v$ be two distinct vertices of $G$  and  let $\sigma \in Aut(G)$. Then $g(u,v)=g(u^{\sigma}, v^{\sigma}).$
\end{lemma}

\begin{proof}
Let $x \in V(G)$. By Lemma \ref{distance preserve}, it follows that $d(x,u)<d(x,v)$ if and only if $d(x^{\sigma}, u^{\sigma})< d(x^{\sigma}, v^{\sigma})$. So, if  $x\in V_{uv}$ then $d(x,u)<d(x,v) \implies d(x^\sigma, u^\sigma)<d(x^\sigma,v^\sigma)\implies x^\sigma \in V_{u^\sigma v^\sigma}$. Thus, we have  $|V_{uv}|\leq |V_{u^\sigma v^\sigma}|$. Similarly, it can be shown that  $|V_{u^\sigma v^\sigma}|\leq |V_{uv}|$. Thus $|V_{uv}|=|V_{u^\sigma v^\sigma}|$ and the result follows.
\end{proof}

\begin{theorem}\label{S_trans}
Let  $G$ be  a connected vertex transitive graph. Then $\mathbb{S}(G)=V(G)$.
\end{theorem}
\begin{proof}
Let $u$ and $v$ be two distinct vertices of  $G$. It is sufficient to show that $s(u)=s(v)$.
Let  $\sigma\in Aut(G)$ such that $u^{\sigma}=v$.  By Lemma \ref{g(u,v)} we have,
\begin{align*}
   s(u)&=\min\{g(u,x): x\in V(G)-u\} \\
         &=\min\{g(u^{\sigma}, x^{\sigma}): x^{\sigma}\in V(G)-u^{\sigma}\}\\                            
         &=\min\{g(v, x^{\sigma}): x^{\sigma}\in V(G)-v\}\\
         &=\min\{g(v, y): y\in V(G)-v\}\\
         &=s(v).
  \end{align*}
 \end{proof}

Is the converse of Theorem \ref{S_trans} true?  To be more precise, we have the following question.
\begin{question}
Let $G$ be a connected graph with $\mathbb{S}(G)=V(G)$. Is $G$  vertex transitive? 
\end{question}

We call the subgraph induced by the vertices of the security center of $G$ as the {\it security subgraph} of $G$. Given a graph $G$, Smart and Slater constructed a graph whose security subgraph is same as $G$.

\begin{proposition} [\cite{Smart}, Theorem 7]
Let $G$ be a graph without isolated vertices. Then there exists a graph $G'$ such that the security subgraph of $G'$ is isomorphic to $G$.
\end{proposition}

In Section \ref{C_center}, we see that the center, median and security center may be different in a graph but all  three coincide for paths and stars. For the path $P_n:12\cdots n$ and the star $K_{1,n-1}$, it is easy to check that 

\begin{equation}\label{CMS_path}
C(P_n)=M(P_n)=\mathbb{S}(P_n)=
\begin{cases}
\{{\frac{n}{2}},{\frac{n}{2}+1}\}  &\mbox{if n is even},\\ 
\{{\frac{n+1}{2}}\}   & \mbox{if n is odd}
\end{cases}
\end{equation}
and 
\begin{equation}\label{CMS_star}
C(K_{1,n-1})=M(K_{1,n-1})=\mathbb{S}(K_{1,n-1})=\{v\}
\end{equation}
\noindent where $v\in V(K_{1,n-1})$ is the vertex of degree $n-1$. 

\section{Three new central parts of a graph}\label{New centers}

In this section, we generalize the definitions of three central parts of  a tree namely the characteristic set, the subtree core and the core vertices. We discuss the analogous centrality properties satisfied by center, median and security center for these new central parts. 

\subsection{The characteristic center}\label{C_center}

Let $G$ be a  graph on $n$ vertices with $V(G)=\{v_1,v_2,\ldots,v_n\}$. The degree matrix $D(G)=(d_{ij})$ of $G$ is the $n\times n$ diagonal matrix with $d_{ii}=deg(v_i)$ for $i=1,2,\ldots,n$. The adjacency matrix $A(G)=(a_{ij})$ of $G$ is the $n\times n$ matrix with $a_{ij}=1$ if $i$ and $j$ are adjacent and $0$ otherwise, for $1\leq i,j \leq n$. The Laplacian matrix $L(G)$ of $G$ is defined as $L(G)=D(G)-A(G).$ It is well known that $L(G)$ is a real symmetric positive semi definite  matrix. The smallest eigenvalue of $L(G)$ is $0$ with all one vector as an eigenvector. The second smallest eigenvalue of $L(G)$  is positive if and only if $G$ is connected (see \cite{Fie1}). The second smallest eigenvalue of $L(G)$ is called the {\it algebraic connectivity} of $G$ and we  denote it by $\mu(G)$. An eigenvector corresponding to $\mu(G)$ is called a {\it Fiedler vector} of $G.$

For a  vertex $v$ of a connected graph $G$,  let  $C_1,C_2,\ldots, C_k$ be the connected components of $G-v$. Note that $k\geq 2$ if and only if $v$ is a cut vertex of $G$. For $i\in \{1,2,\ldots, k\},$ let $\hat{L}(C_i)$ be the principal submatrix of $L(G)$ corresponding to the vertices of $C_i$. Then  $\hat{L}(C_i) $ is invertible and $\hat{L}(C_i)^{-1} $ is a positive matrix (matrix with positive entries) which is called the {\it bottleneck matrix} of $C_i.$ By Perron-Frobenius theorem, $\hat{L}(C_i)^{-1} $ has a simple dominant eigenvalue, called the {\it Perron value} of $C_i$ at $v$. The component $C_j$ is called a {\it Perron component} at $v$ if its Perron value is maximal among the components $C_1, C_2,\ldots, C_k$, at $v$. The next result  describes the entries of  the bottleneck matrices for trees which is useful for our study.

\begin{lemma} [\cite{Kirkland2}, Proposition 1]\label{L-01}
 Let $T$ be a tree and  $v \in V(T).$ Let $T_1$ be a component of $T - v$ and  $L_1$ be the submatrix of $L(T)$ corresponding to $T_1.$ Then $L_1^{-1}=(m_{ij})$, where $m_{ij}$ is the number of edges in common between the paths $P_{iv}$ and $P_{jv},$ where $P_{iv}$ denotes the path joining the vertices $i$ and $v$.
\end{lemma}

Let $Y$ be a Fidler vector of $G$. By $Y(v)$ we mean the co-ordinate of $Y$ corresponding to the vertex $v$ of $G$. A vertex $v$ is called a {\it characteristic vertex}  of $G$ with respect to (w.r.t.) $Y$ if it satisfies one of the following two conditions.
 \begin{enumerate}
 \item[$(i)$] $Y(v)=0$ and there exists a vertex $u$ adjacent to $v$ such that $Y(u)\neq 0$.
 \item[$(ii)$] there exists a vertex $u$ adjacent to $v$ such that $Y(v)Y(u)<0$.
 \end{enumerate}
 \noindent The set of all characteristic vertices of $G$ w.r.t. $Y$ is called the {\it characteristic set} of $G$ w.r.t. $Y$. We denote the characteristic set of $G$ w.r.t.  $Y$ by $\chi(G,Y)$. One of the important reason to consider it as a central part of a tree is the following proposition.

\begin{proposition} [\cite{Fie2}, Theorem 3,14 and \cite{Rm},Theorem 2]\label{Cset}
Let $Y$ be a Fiedler vector of a tree $T$. Then $\chi(T,Y)$ is either a single vertex or two adjacent vertices. Furthermore, $\chi(T,Y)$  is fixed for any Fiedler vector $Y$.
\end{proposition}

Note that for a characteristic vertex $v$, if condition $(ii)$ holds then $u$ is also a characteristic  vertex of $G$ w.r.t. $Y$. In this case, the edge $\{u,v\}$ is known as a {\it characteristic edge} of $G$ w.r.t. $Y$. Thus by a characteristic edge of $G$ w.r.t. $Y$, we mean two adjacent characteristic vertices of $G$ w.r.t. $Y$.  A connection between Perron components and characteristic set of a tree  is described in next three results.
\begin{proposition}[\cite{Kirkland2}, Corollary 1.1]\label{pt1}
Let $T$ be a tree on $n$ vertices. Then the edge $\{i,j\}$ is the characteristic edge of $T$ if and only if the component $T_i$ at vertex $j$ containing the vertex $i$ is the unique Perron  component at $j$ while the component $T_j$ at vertex $i$ containing the vertex $j$ is the unique Perron component at $i.$
\end{proposition}

\begin{proposition}[\cite{Kirkland2}, Corollary 2.1]\label{pt2}
 Let $T$ be a tree on $n$ vertices. Then the vertex $v$ is the characteristic vertex of $T$ if and only if there are two or more Perron components of $T$ at $v$.
\end{proposition}

\begin{proposition} [\cite{Kirkland2}, Proposition 2]\label{pt3} 
Let $T$ be a tree and suppose that $v$ is not a characteristic vertex of $T$. Then the unique Perron component at $v$ contains the characteristic set of $T.$
\end{proposition}

 As an application of Lemma \ref{L-01}, Proposition \ref{pt1} and Proposition \ref{pt2}, we have the following remarks.

\begin{remark}\label{Char_path}
For any Fiedler vector $Y$ of $P_n$, $$\chi(P_n,Y)= \begin{cases}
\{\frac{n}{2},\frac{n}{2}+1\} &\mbox{if $n$ is even},\\
\{\frac{n+1}{2}\} &\mbox{if $n$ is odd}.
\end{cases}$$
 \end{remark} 

\begin{remark}\label{Char_star}
For any Fiedler vector $Y$ of $K_{1,n-1}$,  $\chi(K_{1,n-1},Y)=\{v\}$ where $v$ is the vertex of degree $n-1$ in $K_{1,n-1}.$
\end{remark}

By Proposition \ref{Cset}, it is clear that the characteristic set of a tree is independent of the choice of the Fiedler vector but this is not true for general graphs. We have the following example.

\begin{example}\label{eg_cset}
For the cycle $C_4:v_1v_2v_3v_4v_1$, we have $\mu(C_4)=2.$ It can be easily checked that  $Y_1=(1,0,-1,0)$ and $Y_2=(0,1,0,-1)$  are two Fiedler vectors of $C_4$. So, we get $\chi(C_4,Y_1)=\{v_2,v_4\}$ and $\chi(C_4,Y_2)=\{v_1,v_3\}$.
\end{example}

\noindent The above example shows that, the characteristic set of a graph can be different for different Fiedler vectors. So, claiming it a central part for connected graphs is ill-suited. This motivates us to give a more general definition of the characteristic set to consider it as a central part for graphs.  

\begin{definition}
Let $G$ be a connected graph and $L(G)$ be the Laplacian matrix of $G$. Then the {\bf characteristic center} $\chi(G)$ of $G$ is given by
$$\chi(G)=\{v\in V(G): v \in \chi(G,Y) \;\;\mbox{for some Fiedler vector $Y$}\}.$$
\end{definition}
The term characteristic center (in place of characteristic set) of a tree is first used by Zimmermann in \cite{Zimmermann}. Clearly, $\chi(G)$ is independent  of the choice of the Fiedler vector and for a  tree $T$,  $\chi(T)=\chi(T,Y)$ for any Fiedler vector $Y$. 
So from Remark \ref{Char_path} and Remark \ref{Char_star}, it follows that
$$\chi(P_n)=\begin{cases}
\{\frac{n}{2},\frac{n}{2}+1\} &\mbox{if $n$ is even}, \\
\{\frac{n+1}{2}\} & \mbox{if $n$ is odd}
\end{cases}$$
and $$\chi(K_{1,n-1})=\{v\}$$ where $v$ is the vertex of degree $n-1$ in $K_{1,n-1}$. Also, the following property holds for the characteristic center of a tree.

\begin{proposition} \label{Rmk_Ccen}
The characteristic center of a tree consists of either one vertex or two adjacent vertices.
\end{proposition}

To support the centrality nature of the characteristic center, next we show that $\chi(G)$ lies in a block of $G$.  Let $Y$ be a Fiedler vector of $G$. We call a vertex $v$ has a positive valuation, negative valuation or zero valuation depending upon whether $Y(v)$ is positive, negative or zero, respectively.

\begin{proposition}[\cite{Fie2}, Theorem 3,12]\label{2cases}
Let $G$ be a connected graph and $Y$ be a Fiedler vector of $G$. Then exactly one of the following two cases holds.\\
{\bf Case A:} There is a single block $B_0$ in $G$ which contains vertices with both positive and negative valuations. Each other block contains either only positively valuated  vertices, only negatively valuated vertices or only zero valuated vertices. Every path $P$ starting from $B_0$, which contains at most two cut vertices in each block and exactly one vertex $k$ in $B_0$ has the property that the valuations of the cut vertices of $G$ lying in $P$, form either an increasing or a decreasing  or a zero sequence along this path according to whether $Y(k)>0, Y(k)<0$ or $Y(k)=0$. In the last case all the vertices on $P$ have valuation zero.\\

\noindent {\bf Case B:} No block of $G$ contains both positively and negatively valuated vertices. There exists a unique vertex $z$ of valuation zero which is adjacent to a vertex with non-zero valuation. This vertex $z$ is a cut vertex. Each block contains (with the exception of $z$) either the vertices with positive valuations only, vertices with negative valuations only or vertices with zero valuations only. Every path $P$ starting from $z$ which contains at most two cut vertices in each block has the property that the valuations at its cut vertices either increases and then all valuations of vertices on $P$ are positive(with the exception of $z$), or decreases and then all valuations of the vertices on $P$ are negative (with the exception of $z$) or all valuations of the vertices on $P$ are zero. Every path containing both positively and negatively valuated vertices passes through $z$.
\end{proposition}

 Kirkland and Fallat proved the following result which tells about the position of characteristic center in a graph .
 \begin{lemma}[\cite{Kirkland}, Corollary 2.1]\label{K,F}
Let $G$ be a connected graph. Then either Case A holds for every Fiedler vector, and each such Fiedler vector identifies the same block as being the one  with both positively and negatively valuated vertices, or Case B holds for every Fiedler vector, and each such vector identifies the same vertex $z$ which has zero valuation and is adjacent to one with nonzero valuation.
 \end{lemma}
 
  \noindent The next result follows from Lemma \ref{K,F}.
 \begin{theorem}\label{Block_char}
 The characteristic center of a connected graph $G$ is contained in a block of $G$.
 \end{theorem}
 
We now inspect the characteristic center of connected vertex transitive graphs. Every permutation $\sigma$ on the vertex set $\{v_1,v_2,\ldots,v_n\}$ can be represented by an $n\times n$ permutation matrix $P=(p_{ij})$, where $p_{ij}=1$ if $v_i=\sigma(v_j)$ and $p_{ij}=0$ otherwise. The next lemma is useful for obtaining the characteristic center of a vertex transitive graph and it follows from \cite{Biggs}, Proposition 15.2.
 
 \begin{lemma}\label{PLLP}
 Let $L$ be the Laplacian matrix of $G$ and $\sigma$ be a permutation of $V(G)$. Then $\sigma \in Aut(G)$ if and only if $PL=LP$, where $P$ is the permutation matrix representing $\sigma$.
 \end{lemma}

\begin{theorem}
Let $G$ be a connected vertex transitive graph. Then $\chi(G)=V(G)$.
\end{theorem}
\begin{proof}
Let $u\in \chi(G)$ and $v$ be an arbitrary vertex of $G$. It is sufficient to show that $v\in \chi(G)$. Since $G$ is vertex transitive there exists a $\sigma \in Aut(G)$ such that $u^\sigma=v$. Let $Y$ be a Fiedler vector such that $u\in \chi(G,Y)$ and let $P$ be the permutation matrix representing $\sigma$. Then by Lemma \ref{PLLP},  $PY$ is a Fiedler vector of $G$.  Since $P$ is the matrix representation of $\sigma$, we have $PY(x^\sigma)=Y(x)$ for any $x\in V(G)$.  As $u \in \chi(G)$, we have two possibilities.

\noindent{\bf Case I:} $Y(u)=0$ and there exists $u'$ adjacent to $v$ such that $Y(u')\neq 0$\\
In this case $PY(v)=Y(u)=0$ and $u'^\sigma$ is adjacent to $v$ with $PY(u'^\sigma)=Y(u')\neq 0$.

\noindent{\bf Case II:} There exists $u'$ adjacent to $u$ such that $Y(u)Y(u')< 0$\\
In this case $u'^\sigma$ is adjacent to $v$ and $PY(v)PY(u'^\sigma)=Y(u)Y(u')<0$. 

So, $v\in \chi(G,PY)$ and hence  $v \in \chi(G).$
\end{proof}

Like  the median and security center,  we raise the following question for the characteristic center.
\begin{question}
Let $G$ be a connected graph with $\chi(G)=V(G)$. Is $G$  vertex transitive? 
\end{question}

The following example shows that the characteristic center is different from the   center, median and the security center.

\begin{example}\label{ex1}
The graph $G^*$ in figure \ref{4centers}, has center $C(G^*)=\{4,5,12,13\}$, median $M(G^*)=\{4,13\}$, and security center $\mathbb{S}(G^*)=\{1,4,5,6,7,11,12,13\}$ (See \cite{Sla}, Section 1). It can be checked using Matlab that $\mu(G^*)=0.2$ and $Y=(-0.2485,-0.2837,-0.2623,-0.1883,-0.0166,\\0.1584, 0.3018,0.3772,0.3772,0.3772,0.1584,-0.0166,-0.1883,-0.2623,-0.2837)$ is the unique Fiedler vector up to a scalar multiplication. So it follows that $\chi(G^*)=\{5,6,11,12\}$

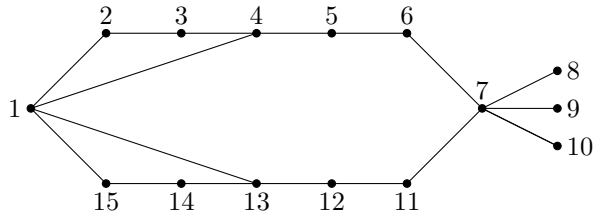
\begin{figure}[h!]
\begin{center}
\begin{tikzpicture}
\filldraw (0,0) node[left]{1} circle [radius=.5mm];
\filldraw (1,1) node[above]{2} circle [radius=.5mm];
\filldraw (1,-1) node[below]{15} circle [radius=.5mm];
\filldraw (2,1) node[above]{3} circle [radius=.5mm];
\filldraw (2,-1) node[below]{14} circle [radius=.5mm];
\filldraw (3,1) node[above]{4} circle [radius=.5mm];
\filldraw (3,-1) node[below]{13}circle [radius=.5mm];
\filldraw (4,1) node[above]{5} circle [radius=.5mm];
\filldraw (4,-1) node[below]{12} circle [radius=.5mm];
\filldraw (5,1) node[above]{6}circle [radius=.5mm];
\filldraw (5,-1) node[below]{11}circle [radius=.5mm];
\filldraw (6,0) node[above]{7} circle [radius=.5mm];
\filldraw (7,0.5) node[right]{8}  circle [radius=.5mm];
\filldraw (7,0) node[right]{9} circle [radius=.5mm];
\filldraw (7,-0.5) node[right]{10}  circle [radius=.5mm];
\draw (7,0.5)--(6,0)--(5,1)--(4,1)--(3,1)--(2,1)--(1,1)--(0,0)--(1,-1)--(2,-1)--(3,-1)--(4,-1)--(5,-1)--(6,0)--(7,-0.5);
\draw (3,-1)--(0,0)--(3,1);
\draw (7,-0.5)--(6,0)--(7,0);
\end{tikzpicture}
\end{center}
\caption{The graph $G^*$ with distinct central parts}\label{4centers}
\end{figure}
\end{example}

We call the subgraph induced by $\chi(G)$, the characteristic subgraph of $G$. We denote the set of all $n\times n$ complex matrices by $M_n$. Next we prove that for any graph $G$ there exists a super graph with $G$ as a characteristic subgraph.  In that context the following two lemmas are useful.

\begin{lemma}[\cite{Horn}, Corollary 4.3.3]\label{Weyl}
Let $A, B \in M_n$ be Hermitian matrices. Let  $\lambda_1(A+B)\leq \lambda_2(A+b)\leq \cdots \leq  \lambda_n(A+B)$ be the eigenvalues of $A+B$ and $\lambda_1(A)\leq \lambda_2(A)\leq \cdots \leq \lambda_n(A)$ be the eigenvalues of $A$. Suppose that $B$ has no negative eigenvalues. Then $\lambda_i(A+B)\geq \lambda_i(A)$ for $i=1,2,\ldots,n$.
\end{lemma}

\begin{corollary}\label{CWeyl}
Let $u$ and $v$ be two nonadjacent vertices of $G$ and let $\lambda(G)\leq \lambda_2(G)\leq \cdots \leq \lambda_n(G)$ be the Laplacian eigenvalues of $G$. Suppose $G'$ is the graph obtained from $G$ by joining $u$ and $v$ with an edge. Then $\lambda_i(G)\leq \lambda_i(G')$ for $i=1,2,\ldots,n$.
\end{corollary}

\begin{lemma}[\cite{Horn}, Theorem 4.3.28]  \label{CIT}
Supose $A\in M_n$ is a Hermitian matrix. Let $B \in M_m$ with $m<n$ be a principal submatrix of $A$. Suppose $A$ has eigenvalues  $\alpha_1\leq \alpha_2\leq \cdots \leq \alpha_{n}$  and $B$ has eigenvalues $\beta_1\leq\beta_2\leq \cdots \leq \beta_{m}$. Then $\alpha_i \leq \beta_i \leq \alpha_{i+n-m}$ for $i=1,2,\ldots m$.  
\end{lemma}

\begin{theorem}
For any graph $G$ (may be disconnected) there exist a graph $G_{ch}$ such that the characteristic subgraph of $G_{ch}$ is isomorphic to $G$.
\end{theorem}
\begin{proof}
Let $G$ be a graph with $V(G)=\{v_1,v_2,\ldots, v_n\}$. Consider the graph $G_{ch}$ with $V(G_{ch})=V(G)\cup \{u_1,u_2,u_3,u_4\}$ and edges of $G_{ch}$ are the edges of $G$ together with $2n+2$ new edges $\{u_1,u_2\},\{u_3,u_4\}, \{u_2,v_i\}$ and $\{u_3,v_i\}$, $i=1,2\ldots, n$. Note that for any $G$, $G_{ch}$ is connected.

\begin{figure}[h!]
\begin{center}
\begin{tikzpicture}
\draw (0,0) circle [radius=1cm];
\filldraw (-3.5,0) node[above]{$u_1$} circle [radius=.5mm] (-2.5,0) node[above]{$u_2$} circle [radius=.5mm] (2.5,0) node [above]{$u_3$} circle [radius=.5mm] (3.5,0) node[above]{$u_4$} circle [radius=.5mm];
\draw (-3.5,0)--(-2.5,0)--(-.5,.867) (.5,.867)--(2.5,0)--(3.5,0) (-2.5,0)--(-.5,-.867)  (.5,-.867)--(2.5,0);
\draw (0,0) node {$G$}  (0,-1.5) node{$G_{ch}$};
\end{tikzpicture}
\caption {The graph $G_{ch}$}
\end{center}
\end{figure}
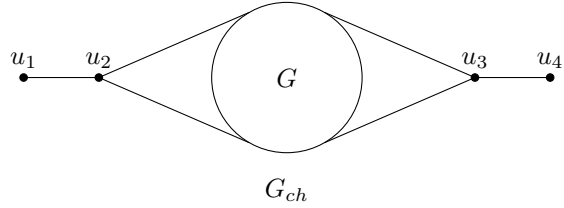

Taking the vertices of $G_{ch}$ in the order $u_1,u_2, v_1,\ldots, v_n, u_3,u_4$, the Laplacian matrix of $G_{ch}$ is \\

\begin{equation*}
L(G_{ch})=\left(
\begin{array}{ccc}
      \begin{matrix}
       1&-1\\
      -1&n+1
      \end{matrix} & \begin{matrix}
                            0&\cdots &0\\
                           -1&\cdots &-1
                           \end{matrix} & \begin{matrix}
                                                  0 &0\\
                                                  0 &0
                                                 \end{matrix}\\

\begin{matrix}
   0&-1\\
   \vdots& \vdots\\
    0&-1
\end{matrix}& \begin{pmatrix}
                        &&\\
                        &L(G)+2I&\\
                         &&
                       \end{pmatrix}&  \begin{matrix}
                                                 -1&0\\
                                            \vdots&\vdots\\
                                                  -1&0
                                               \end{matrix}\\
\begin{matrix}
0&0\\
0&0
\end{matrix}& \begin{matrix}
                       -1&\cdots&-1\\
                        0&\cdots &0
                       \end{matrix}& \begin{matrix}
                                              n+1&-1\\
                                             -1&1
                                             \end{matrix}

\end{array}
\right).
\end{equation*}
Let $\lambda_1\leq \lambda_2 \leq \ldots \leq \lambda_{n+4}$ be the eigenvalues of $L(G_{ch})$. It can be checked that $\lambda=\frac{n+2-\sqrt{n^2+4}}{2}$ is an eigenvalue of $L(G_{ch})$ with corresponding eigenvector $Y=(1, 1-\lambda , 0, \cdots ,0, \lambda-1, -1)^T$. As $n\geq 1$, $0<\lambda<1$. Now we show that $\lambda_2=\lambda$ and multiplicity of $\lambda_2$ is one.

First take  $G=\bar{K_n}$, where $\bar{K_n}$ is the complement of the complete graph on $n$ vertices. Let $\alpha_1\leq \alpha_2\leq \cdots \leq \alpha_{n+4}$ be the eigenvalues of $L(\bar{K_n}_{ch})$. Note that $\lambda$ is an eigenvalue of $L(\bar{K_n}_{ch})$ with corresponding eigenvector $Y=(1,1-\lambda , 0, \cdots ,0, \lambda-1, -1)^T$. Consider the principal submatrix $L(\bar{K_n}_{ch})(2,n+3)$ of $L(\bar{K_n}_{ch})$, obtained by deleting $2^{nd}$ and $(n+3)^{th}$ rows and columns, respectively. Then $L(\bar{K_n}_{ch})(2,n+3)$ is the diagonal matrix $diag(1,2,2,\ldots,2,1).$ Let  $\beta_1\leq\beta_2\leq \cdots \leq \beta_{n+2}$ be the eigenvalues of $L(\bar{K_n}_{ch})(2,n+3)$.   Taking $i=1$ in Lemma \ref{CIT}, we get $\alpha_3 \geq \beta_1=1$. Since $\alpha_1=0$, we have $\alpha_2=\lambda$.

By Corollary \ref{CWeyl}, $\lambda_i\geq \alpha_i$ for $i=1,2,\ldots,n+4$. So, $\lambda_3\geq 1$ and hence $\lambda_2=\lambda$ is an  eigenvalue of $L(G_{ch})$ with multiplicity one. Therefore, $Y$ is a Fiedler vector of $G_{ch}$ and we get $\chi(G_{ch},Y)=V(G)$. Since multiplicity of $\lambda$ is one,  so every Fiedler vector of $G_{ch}$ is a scalar multiple of $Y$. Hence, $\chi(G_{ch})=\chi(G_{ch},Y)=V(G)$ and the result follows.

\end{proof}

\subsection{The subgraph core} \label{S_core}

In $2005$, S\'zekely and Wang \cite{Sze} defined a new central part of a tree different from both the center and the centroid. For a vertex $v$  of a tree $T$, the {\it subtree number} $f_T(v)$ of $v$ is the number of subtrees of $T$ containing $v$. The {\it subtree core} of $T$ is the set of vertices having maximum subtree number. The following result is a motivation towards considering the subtree core as a  central part of a tree.

\begin{proposition}[\cite{Sze}, Theorem 9.1] \label{prop:Sc}
The subtree core of a tree consists of either a single vertex or two adjacent vertices.
\end{proposition}

The subtree core is exclusively defined for trees. We give a very natural extension of the subtree core of a tree  to general graphs. As trees have subtrees, graphs have connected subgraphs. So we define the subgraph core of a graph as follow. 

\begin{definition}
Let $G$ be a graph and $v\in V(G).$ The subgraph number $f_G(v)$ of $v$ is the number of connected subgraphs of $G$ containing $v.$ The set of vertices of $G$  having maximum subgraph number is called the {\bf subgraph core} of $G$ and we denote it by $S_c(G)$.
\end{definition}

\noindent Note that for a tree $T$, the subgraph core of $T$ is same as the subtree core of $T$. So the following basic property for a central part in trees hold.

\begin{proposition}
The subgraph core of a tree consists of either a single vertex or two adjacent vertices.
\end{proposition}  
 Also, we have

$$S_c(P_n)=\begin{cases}
\{\frac{n}{2},\frac{n}{2}+1\} &\mbox{if $n$ is even}, \\
\{\frac{n+1}{2}\} & \mbox{if $n$ is odd}
\end{cases}$$
and $$S_c(K_{1,n-1})=\{v\}$$ where $v$ is the vertex of degree $n-1$ in $K_{1,n-1}$.

An important property of a central part of a graph is that it lies in a block. We strongly feel that this is true for the subgraph core also. So, we conjecture the following.

\begin{conjecture}
The subgraph core of a graph $G$ is contained in a block of $G$.
\end{conjecture}

In the next result we justify the centrality property which characterises the subgraph core of vertex transitive graphs.

\begin{theorem}
Let $G$ be a connected vertex transitive graph. Then $S_c(G)=V(G)$.
\end{theorem}

\begin{proof}
Let $u,v\in V(G)$. It is sufficient to show that $f_G(u)=f_G(v)$. Since $G$ is vertex transitive, there exists $\sigma\in Aut(G)$ such that $u^{\sigma}=v$. Let $S_u$ be the set of all connected subgraphs of $G$ containing $u$ and $S_v$ be the set of all connected subgraphs of $G$ containing $v$. Then there is a bijection between $S_u$ and $S_v$ which sends $H$ to $H^\sigma$ where $H^\sigma$ is the subgraph of $G$ with $V(H^\sigma)=\{x^{\sigma}: x\in V(H)\}$ and $E(H^\sigma)=\lbrace \{ x^{\sigma}, y^{\sigma} \}: \{ x,y \} \in E(H) \rbrace$. This implies $|S_u|=|S_v|$ and hence, $f_G(u)=f_G(v)$.
\end{proof}

We have  the following questions for the subgraph core.
  \begin{question}
Let $G$ be a connected graph with $S_c(G)=V(G)$. Is $G$ vertex transitive? 
  \end{question}
We call the subgraph induced by the subgraph core of $G$ as the {\it core subgraph} of $G$. The following question can be asked regarding the core subgraph.  
 \begin{question}
  Given a graph $G$ does there exist a graph $G'$ such that the core subgraph of $G'$ is isomorphic to $G$?
 \end{question}
 
 \subsection{The core vertices}
 
The core vertices of a tree is the most recently defined central part introduced by Zhang et al. in \cite{Zhang}. Let $T$ be a tree and $v\in V(T)$. The {\it eccentric subtree number} $\epsilon(v)$ is defined as $\epsilon(v)=\epsilon_T(v)=\min\{f_T(v,u): u\in V(T)\}$ where $f_T(v,u)$ denotes the number of subtrees of $T$ containing both $v$ and $u$. The {\it core vertices} of $T$ is the set of vertices having maximum eccentric subtree number.

 The core vertices is defined for trees only but it has a natural extension to connected graphs. So we define the core vertices of a graph as follow.
 
 \begin{definition}
 Let $G$ be a graph. We define the {\it eccentric subgraph number} $\epsilon(v)$ of $v$ in $G$ as $\epsilon(v)=\epsilon_G(v)=\min\{f_G(v,u): u\in V(G)\}$ where $f_G(v,u)$ denotes the number of connected subgraphs of $G$ containing both $v$ and $u$. The {\bf core vertices} of $G$ is the set of vertices having maximum eccentric subgraph number. We denote the core vertices of $G$ by $\mathcal{C}(G) $.
 \end{definition}

\begin{proposition}[\cite{Zhang}, Theorem 3.4]
The core vertices of a tree consists of either one vertex or two adjacent vertices.
\end{proposition}
 
By $f_T(x,y,z)$ we denote the number of subtrees of $T$ containing the vertices $x,y$ and $z$ and by $f_T(x,y,\bar{z})$, the number of subtrees of $T$ containing $x$ and $y$ but not $z$.

\begin{proposition}\label{minf(u,v)}
Let $v\in V(T)$ and $\min\{f_T(v,x):x \in V(T)\}=f_T(v,u)$. Then $u$ must be a pendant vertex.
\end{proposition}

\begin{proof}
Suppose $u$ is a non pendant vertex. Then there exists a pendant vertex $w$ such that $u$ lies on the path joining $v$ and $w$.  Now
 $f_T(v,u)=f_T(v,u,w)+ f_T(v,u,\bar{w})=f_T(v,w)+f_T(v,u,\bar{w})>f_T(v,w)$, which is a contradiction.
\end{proof}
 
 For the path $P_n:12\cdots n$, we have  $  \epsilon(k)=f_{P_n}(k,n)=k=\epsilon(n-k+1) \;\; \mbox{for} \;\; 1\leq k\leq \left\lceil\frac{n}{2}\right\rceil$ and for the star  $K_{1,n-1}$,  $n\geq 3$ we have, $\epsilon(v)=2^{n-2}$ where $v$ is the vertex of degree $n-1$ and $\epsilon(u)=2^{n-3}$ for any other vertex $u$. So, it follows that  
 
 $$\mathcal{C}(P_n)=\begin{cases}
\{\frac{n}{2},\frac{n}{2}+1\} &\mbox{if $n$ is even}, \\
\{\frac{n+1}{2}\} & \mbox{if $n$ is odd}
\end{cases}$$
and $$\mathcal{C}(K_{1,n-1})=\{v\}$$ where $v$ is the vertex of degree $n-1$.

\begin{proposition}[\cite{Zhang}, Theorem 3.1]\label{Concave}
Suppose $u,v,w\in V(T)$ such that $\{u,v\}, \{v,w\}\in E(T)$. Then $2\epsilon(v)\geq \epsilon(u)+\epsilon(w)$ with a possible equality only if $deg(v)=2$. 
\end{proposition}

 Now we give an example of a tree in which the center, centroid, characteristic center, subgraph core and core vertices are disjoint.
 
  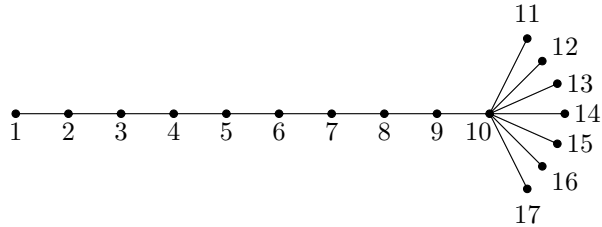
\begin{figure}[h] 
 	\begin{center}
 		\begin{tikzpicture}
 		\filldraw (0,0)node [below]{1} circle [radius= .5mm]--(.7,0)node [below]{2} circle[radius=.5mm]--(1.4,0)node [below]{3} circle [radius= .5mm]--(2.1,0)node [below]{4} circle[radius=.5 mm]--(2.8,0)node [below]{5} circle [radius= .5mm]--(3.5,0)node [below]{6} circle[radius=.5 mm]--(4.2,0) node [below]{7}circle [radius= .5mm]--(4.9,0)node [below]{8} circle[radius=.5 mm]--(5.6,0) node [below]{9}circle [radius= .5mm]--(6.3,0)circle[radius=.5 mm]--(7.3,0)node [right]{14} circle[radius=.5 mm];
 			\draw (6.15,0) node [below]{10};
 			\filldraw (6.8,1)node[above,outer sep=3pt]{11}circle [radius=.5mm]--(6.3,0);
 			\filldraw (7,.7) circle [radius=.5mm]--(6.3,0);
 			\filldraw (7.2,.4) node[right]{13}circle [radius=.5mm]--(6.3,0);
 			\draw (7,.9)node[right]{12};
 			\filldraw (6.8,-1)node[below,outer sep=3pt]{17}circle [radius=.5mm]--(6.3,0);
 			\filldraw (7,-.7) circle [radius=.5mm]--(6.3,0);
 			\filldraw (7.2,-.4) node[right]{15}circle [radius=.5mm]--(6.3,0);
 			\draw (7,-.9)node[right]{16};
 			\end{tikzpicture}
 		\end{center}
		\vskip -0.7cm
 		\caption{ A tree $T$ with disjoint central parts}\label{P10,7}
 	\end{figure}

 \begin{example}\label{ex2}
For the tree $T$ in Figure \ref{P10,7}, observe that $f_T(3,11)=f_T(3,i)$ for $i=12,13,\ldots,17$. So, from Proposition \ref{minf(u,v)}, it is clear that $\min\{f_T(3,v): v\in V(T)\}= \min\{f_T(3,1),f_T(3,11)\}$. We have $f_T(3,1)=f_T(3,1,\bar{10})+f_T(3,1,10)=7+2^7=135$ and 
 $f_T(3,11)=f_T(3,11,\bar{2})+f_T(3,11,2)=2^6+2\times 2^6=3\times2^6=192.$ Therefore, $\epsilon(3)=135$. Similar calculations give $\epsilon(2)=128$ and $\epsilon(4)=134.$ So, by Proposition \ref{Concave}, it follows that $\mathcal{C}(T)=\{3\}$. Also we have (see \cite{Pandey2}, Example 1.6) $C(T)=\{6\}$, $C_d(T)=\{9\}$, $S_c(T)=\{10\}$ and $\chi(T)=\{7,8\}$. 
 \end{example}
Thus from Example \ref{ex1} and Example \ref{ex2}, it can be  observed that the center, median, security center, characteristic center, subgraph core and the core vertices may be all different in  a graph. The following result justifies the centrality nature of the core vertices for connected graphs

\begin{theorem}
Let $G$ be a connected vertex transitive graph. Then $\mathcal{C}(G)=V(G)$.
\end{theorem}
\begin{proof}
Let $u,v\in V(G)$. It is sufficient to show that $\epsilon(u)=\epsilon(v)$. Let $S_{xy}$ be the set of all connected subgraphs of $G$ containing both the vertices $x$ and $y$. Then for any $\tau\in Aut(G)$, The map $\phi : S_{xy}\rightarrow S_{x^\tau y^\tau}$ defined as $\phi(H)=H^\tau$ is a bijection, where $H^\tau$ is the graph having $V(H^\tau)=\{a^\tau: a\in V(H)\}$ and $E(H^\tau)=\lbrace \{a^\tau,b^\tau\}: \{ a,b \} \in E(H) \rbrace$.

Since $G$ is vertex transitive there exists $\sigma\in Aut (G)$ such that $u^\sigma=v$. Now 
\begin{align*}
\epsilon(u)&= \min\{f(u,x): x\in V(G)\}\\ 
                 &=\min\{ |S_{ux}|: x\in V(G)\}\\
                 &=\min\{|S_{u^\sigma x^\sigma}|: x\in V(G)\}\\
                 &=\min\{|S_{v x^\sigma}|: x\in V(G)\}\\
                 &=\min\{|S_{vy}|: y\in V(G)\}\\
                 &=\min\{f(v,y):y\in V(G)\}\\
                 &=\epsilon(v).
\end{align*}
\end{proof}
 We end this section with the following questions related to the core vertices of a graph.
  \begin{question}
  Let $G$ be a connected graph with $\mathcal{C}(G)=V(G)$. Is $G$ vertex transitive? 
  \end{question}
\begin{question}
Is the core vertices of a connected graph $G$ contained in a block of $G$?
\end{question}
 \begin{question}
 \it Given a graph $G$ does there exists a graph $G'$ such that the subgraph induced by the core vertices of $G'$ is isomorphic to $G$?
 \end{question}
 
\section{Conclusion}\label{FQ}

We have discussed many properties of the center, median, security center, characteristic center, subgraph core and the core vertices of a graph and shown that all six may be different in a graph. It is observed that all these central parts have similar behaviour in a graph.  We conclude that  following are the properties of any central part of a connected graph.

\begin {itemize}
\item a central part of a tree is either a single vertex or two adjacent vertices.
\item any central part of path and star coincide with the center of the path and star respectively.
\item a central part of a graph always lie in a block.
\item a central part of any vertex transitive graph consists of all the vertices of the graph. 
\end{itemize}

\end{document}